\documentclass[12pt]{article}

\newtheorem{thm}{Theorem}
\newtheorem{alg}[thm]{Algorithm}

\newtheorem{cnj}[thm]{Conjecture}
\newtheorem{cor}[thm]{Corollary}

\newtheorem{lem}[thm]{Lemma}
\newtheorem{prb}[thm]{Problem}

\usepackage{graphicx,latexsym,amsmath,amsfonts,amssymb,mathrsfs}

\def\a{{\alpha}}
\def\b{{\beta}}
\def\c{{\gamma}}
\def\d{{\delta}}

\def\k{{\kappa}}
\def\l{{\lambda}}

\def\w{{\omega}}

\def\bT{{\mathbf T}}

\def\zN{{\mathbb N}}

\def\cF{{\cal F}}

\def\sNP{{\sf NP}}

\def\sw{{\sf w}}

\def\Cp{{C^\pr}}
\def\crr{{\ct^r}}
\def\crz{{\ct_0}}
\def\ct{{\rm x}}
\def\ctp{{\ct^\pr}}
\def\dpr{{\d^\pr}}

\def\ds{{\d^*}}

\def\dtr{{D_3(r)}}
\def\Gp{{G^\pr}}
\def\Kp{{K^\pr}}

\def\Nx!{{N_x^!}}

\def\rp{{r^\pr}}
\def\sp{{s^\pr}}
\def\SP{{S^+}}

\def\up{{u^\pr}}
\def\vp{{v^\pr}}
\def\wp{{w^\pr}}

\def\Xrz{{X_0}}

\def\Ars{{A_{rs}}}
\def\crs{{\ct_{rs}}}
\def\Crs{{C_{rs}}}
\def\Krs{{K_{rs}}}
\def\Srs{{S_{rs}}}
\def\trs{{t_{rs}}}
\def\trsp{{t_{\rp\sp}}}
\def\Xrs{{X_{rs}}}

\def\deg{{\sf deg}}
\def\diam{{\sf diam}}
\def\dist{{\sf dist}}

\def\ecc{{\sf ecc}}

\def\mt{{\emptyset}}

\def\per{{\psi}}
\def\perp{{\psi^\pr}}
\def\phx{{\phi}}
\def\ppr{{\pr\pr}}
\def\pr{{\prime}}
\def\Pf{{\it Proof.\ \ }}
\def\pf{{\hfill $\Box$\bigskip}}

\def\rar{{\rightarrow}}

\def\sse{{\subseteq}}

\def\sqr#1#2{{\vcenter{\hrule height.#2pt
        \hbox{\vrule width.#2pt height#1pt \kern#1pt
                \vrule width.#2pt}
        \hrule height.#2pt}}}

\def\str{{$\star$}}
\def\sstr{{$\star\star$}}

\def\lf{{\lfloor}}
\def\rf{{\rfloor}}

\begin{document}

%
%
\title{Pebbling in Split Graphs\\
\ \\}

\author{
Liliana Alc\' on\thanks{
Departamento de Matem\' atica,
Facultad Ciencias Exactas,
Universidad Nacional de La Plata, 
CC 172, (1900) La Plata
Argentina
}\\
Marisa Gutierrez\footnotemark[1]
\ \thanks{
CONICET, Argentina
}\\
Glenn Hurlbert\thanks{
School of Mathematics and Statistics,
Arizona State University,
Tempe, AZ 85287 USA
}
\thanks{
Research supported by a J. William Fulbright Scholarship while
on sabbatical at UNLP
}
}
\maketitle

\newpage

%
%
\begin{abstract}

Graph pebbling is a network optimization model for transporting
discrete resources that are consumed in transit: the movement of
two pebbles across an edge consumes one of the pebbles.  The pebbling
number of a graph is the fewest number of pebbles $t$ so that, from
any initial configuration of $t$ pebbles on its vertices, one can
place a pebble on any given target vertex via such pebbling steps.
It is known that deciding if a given configuration on a particular
graph can reach a specified target is {\sf NP}-complete, even for
diameter two graphs, and that deciding if the pebbling number has a
prescribed upper bound is $\Pi_2^{\sf P}$-complete.

On the other hand, for many families of graphs there are formulas
or polynomial algorithms for computing pebbling numbers; for example,
complete graphs, products of paths (including cubes), trees, cycles,
diameter two graphs, and more.  Moreover, graphs having minimum
pebbling number are called Class 0, and many authors have studied
which graphs are Class 0 and what graph properties guarantee it,
with no characterization in sight.

In this paper we investigate an important family of diameter three
chordal graphs called split graphs; graphs whose vertex set can be
partitioned into a clique and an independent set.  We provide
a formula for the pebbling number of a split graph, along with an
algorithm for calculating it that runs in $O(n^\b)$ time, where
$\b=2\w/(\w+1)\cong 1.41$ and $\w\cong 2.376$ is the exponent
of matrix multiplication.  Furthermore we determine that all
split graphs with minimum degree at least 3 are Class 0.
\medskip

\smallskip


\end{abstract}

\newpage

%
%
\section{Introduction}\label{Intro}

Graph pebbling is a network optimization model for transporting
discrete resources that are consumed in transit: while two pebbles
cross an edge of a graph, only one arrives at the other end as the
other is consumed (or lost to a toll, one can imagine).  This
operation is called a {\it pebbling step}.  The basic questions in
the subject revolve around deciding if a particular configuration
of pebbles on the vertices of a graph can {\it reach} a given {\it
root} vertex via pebbling steps (for this reason, graph pebbling
is carried out on connected graphs only).  If a configuration can
reach $r$, it is called $r$-{\it solvable}, and $r$-{\it unsolvable}
otherwise.

Various rules for pebbling steps have been studied for years and
have found applications in a wide array of areas.  One version,
dubbed {\it black and white} pebbling, was applied to computational
complexity theory in studying time-space tradeoffs (see
\cite{HoPaVa,PatHew}), as well as to optimal register allocation
for compilers (see \cite{Seth}).  Connections have been made also
to pursuit and evasion games and graph searching (see \cite{KirPip,Pars}).
Another ({\it black} pebbling) is used to reorder large sparse
matrices to minimize in-core storage during an out-of-core Cholesky
factorization scheme (see \cite{GiLeTa,Klaw,Liu}).  A third version
yields results in computational geometry in the rigidity of graphs,
matroids, and other structures (see \cite{GurShe,StrThe}).  The
rule we study here originally produced results in combinatorial
number theory and combinatorial group theory (the existence of zero
sum subsequences --- see \cite{Chun,EllHur}) and have recently been
applied to finding solutions in $p$-adic diophantine equations (see
\cite{Knap}).  Most of these rules give rise to computationally
difficult problems, which we discuss for our case below.

We follow fairly standard graph terminology (e.g. \cite{West}),
with a graph $G=(V,E)$ having $n=n(G)$ vertices $V=V(G)$ and having
edges $E=E(G)$.  The {\it eccentricity} $\ecc(G,r)$ for a vertex
$r\in V$ equals $\max_{v\in V}\dist(v,r)$, where $\dist(x,y)$ denotes
the length (number of edges) of the shortest path from $x$ to $y$;
the {\it diameter} $\diam(G)=\max_{r\in V}\ecc(G,r)$.  When $G$ is
understood we will shorten our notation to $\ecc(r)$.

The most studied graph pebbling parameter, and the one investigated
here, is the {\it pebbling number} $\pi(G)= \max_{r\in V}\pi(G,r)$,
where $\pi(G,r)$ is defined to be the minimum number $t$ so that
every configuration of size at least $t$ is $r$-solvable.  The {\it
size} $|C|$ of a configuration $C:V\rar\zN=\{0,1,\ldots\}$ is its
total number of pebbles $\sum_{v\in V}C(v)$.  Simple lower bounds
like $\pi(G)\ge n$ (sharp for complete graphs, cubes, and,
probabilisticaly, almost all graphs) and $\pi(G)\ge 2^{\diam(G)}$
(sharp for paths and cubes, among others) are easily derived.  Graphs
satisfying $\pi(G)=n$ are called {\it Class $0$} and are a topic
of much interest (e.g. \cite{BekCus, BlaSch, ClHoHu, CHHM, CzyHur,
CHKT}).  Surveys on the topic can be found in \cite{HurlRPGP, HurlGGP,
HurlPP}, and include variations on the theme such as $k$-pebbling,
fractional pebbling, optimal pebbling, cover pebbling, and pebbling
thresholds, as well as applications to combinatorial number theory
and combinatorial group theory (see references).

Computing graph pebbling numbers is difficult in general.  The
problem of deciding if a given configuration of pebbles on a graph
can reach a particular vertex was shown in \cite{HeHeHu,HurKie} to
be \sNP-complete (via reduction from the problem of finding a perfect
matching in a $4$-uniform hypergraph).  The problem of deciding if
a graph $G$ has pebbling number at most $k$ was shown in \cite{HeHeHu}
to be $\Pi_2^{\sf P}$-complete.\footnote{That is, complete for the
class of problems computable in polynomial time by a {\sf co-NP}
machine equipped with an oracle for an {\sf NP}-complete language.}

On the other hand, pebbling numbers of many graphs are known: for
example, cliques, trees, cycles, cubes, diameter two graphs, graphs
of connectivity exponential in its diameter, and others.  In
particular, in \cite{PaSnVo} the pebbling number of a diameter $2$
graph $G$ was determined to be $n$ or $n+1$.  Moreover, \cite{ClHoHu}
characterized those graphs having $\pi(G)=n+1$ (a slight error in
the characterization was corrected by \cite{BlaSch}).  All such
connectivity $1$ graphs have $\pi(G)=n+1$.  The smallest such
$2$-connected graph is the {\it near-Pyramid} on $6$ vertices, which
is the $6$-cycle $(r,a,p,c,q,b)$ with an extra two or three of the
edges of the triangle $(a,b,c)$ (the {\it Pyramid} has all three).
All diameter $2$ graphs with pebbling number $n+1$ can be described
by adding simple structures to the near-Pyramid.  It was shown in
\cite{HeHeHu} that one can recognize such graphs in quartic time.

Here we begin to study for which graphs their pebbling numbers can
be calculated in polynomial time.  Aiming for tree-like structures
(as was considered in \cite{CHHM}), one might consider chordal
graphs of various sorts.  Moving away from diameter $2$, one might
consider diameter $3$ graphs; recently (\cite{PoStYe}), the tight
upper bound of $\lf 3n/2\rf+2$ has been shown for this class.
Combining these two thoughts we study split graphs in this paper,
and find that their pebbling numbers can be calculated quickly, in
fact, in $O(n^{1.41})$ time.\footnote{Here $\b\cong 1.41$ satisfies
$\b=2\w/(\w+1)$, where $\w\cong 2.376$ is the exponent of matrix
multiplication.}

Split graphs can be described by adding simplicial vertices ({\it
cones}) to a fixed clique.  In other words, a graph is a {\it split}
graph if its vertices can be partitioned into an independent set
$S$ and a clique $K$.  Notice that the Pyramid is a split graph
with clique $\{a,b,c\}$ and cone vertices $r$, $p$, and $q$.  The
Pyramid plays a key role in the theory of split graphs.  However,
the Pyramid has diameter $2$, and we are interested in diameter $3$
split graphs.

\begin{figure}
\centerline{\includegraphics[height=1.8in]{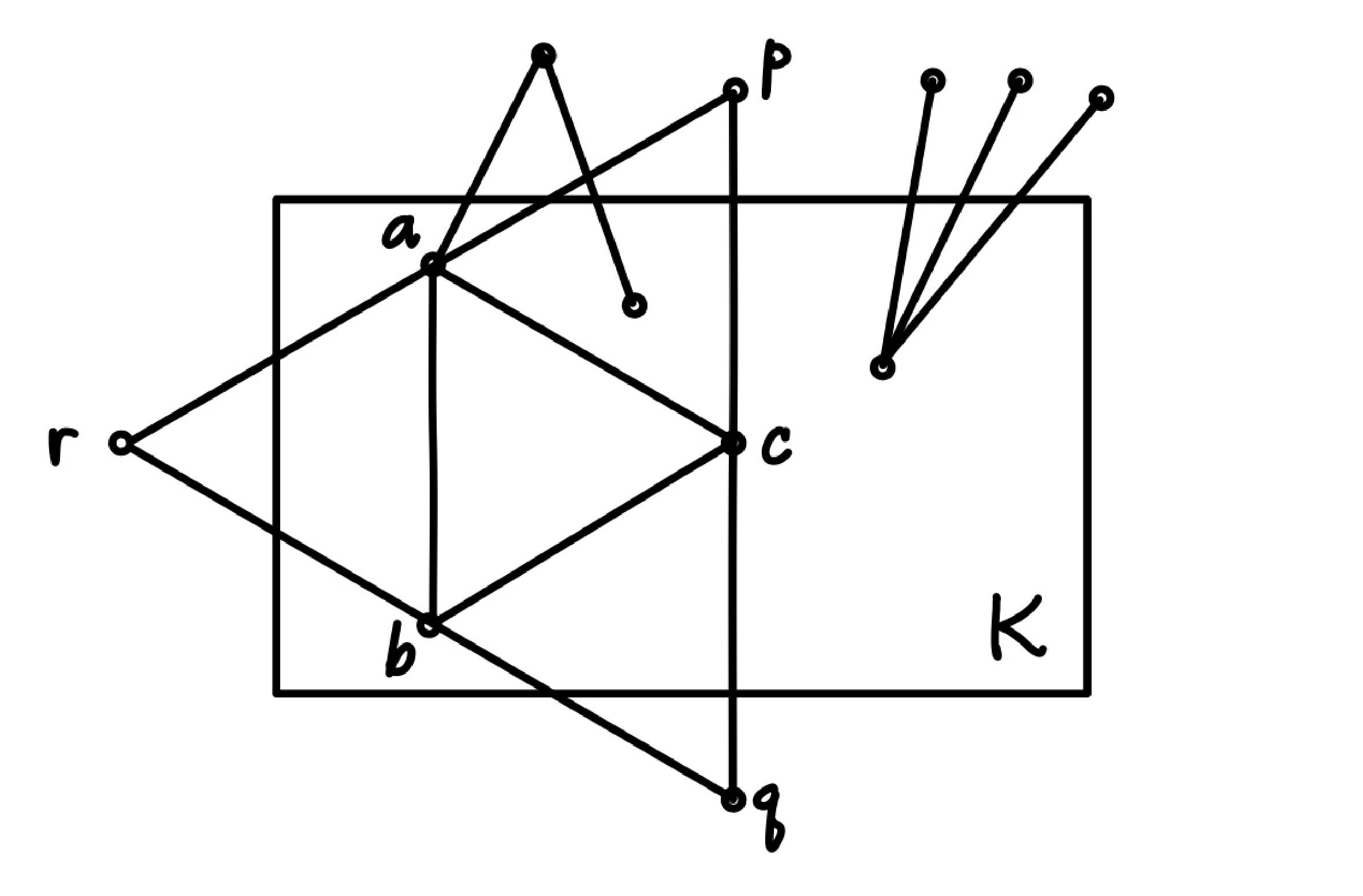}$\quad$
\includegraphics[height=1.8in]{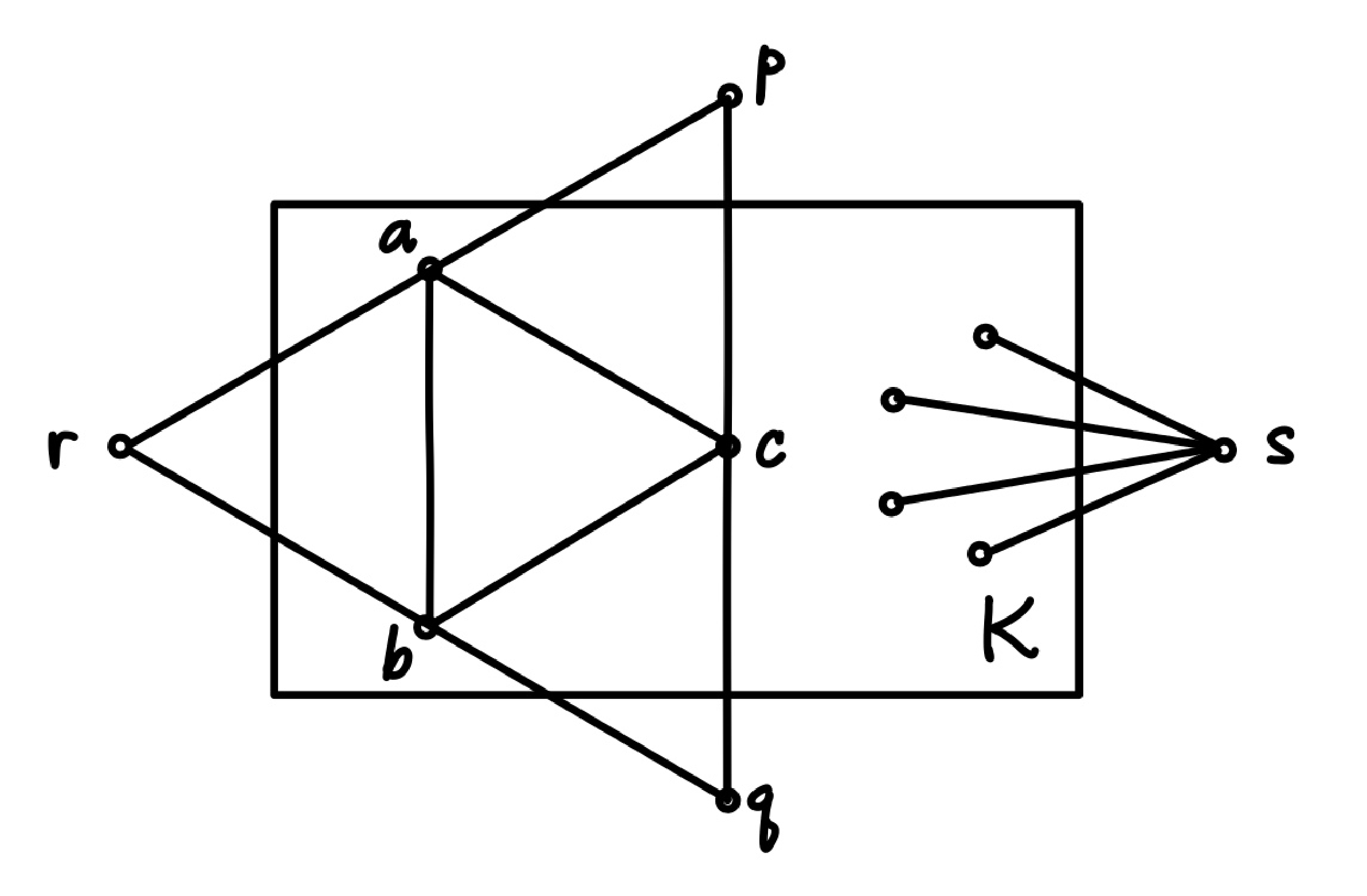}}
\caption{Examples of Pereyra (left) and Phoenix (right) graphs}\label{PerPho}
\end{figure}

It turns out that Pereyra and Phoenix graphs (which we define below
and necessarily contain the Pyramid) are important for our work
(see Fig. \ref{PerPho}).  We say that $G$ has a Pyramid if there
exist three cone vertices with degree $2$ whose neighborhoods do
not have the Helly property (that is, their neighborhoods form a
triangle).  We say that the subgraph induced by the closed neighborhoods
of the three cone vertices is a Pyramid of $G$.  If $r$ is one of
the three cone vertices we say it is an $r$-Pyramid.  A graph $G$
is called $r$-{\it Pereyra} if it has an $r$-Pyramid, none of whose
vertices is a cut vertex of $G$.  Denote by $\ds(G,r)$ the minimum
degree among all vertices at maximum distance from $r$.  A graph
$G$ is $r$-{\it Phoenix} if it is $r$-Pereyra, $\ecc(r)=3$, and
$\ds(G,r)\ge 4$.  A {\it Pereyra} (resp. {\it Phoenix}) graph is
$r$-Pereyra (resp. $r$-Phoenix) for some $r$.

Like the Pyramid, an $r$-Pereyra graph having $\ecc(r)=2$ has
pebbling number one more than ``normal''; that is, it is an exception
to how most of the graphs in its class behave.  On such $G$, the
configuration that places $3$ pebbles on $p$ and $q$, $0$ pebbles
on $r$, $a$, $b$, and $c$, and $1$ everywhere else is $r$-unsolvable,
showing that $\pi(G,r) \ge n + \ct + 1$, where $\ct$ is the number
of cut vertices of $G$.  (In the course of proving Theorem \ref{rootE2}
below, one finds that this configuration is the unique $r$-unsolvable
configuration of size $n + \ct$ on $G$.) We will find analogous
behavior for $r$-Phoenix graphs as well.


Notationally, we abbreviate $\deg(x)$ by $d_x$.  We also abbreviate
$N(x)$ by $N_x$ (so that $d_x=|N_x|$), with $[N_x]$ denoting
$N_x\cup\{x\}$.  If $v\in S$, we define $K_v=K-N_v$.  We denote the
set of cut vertices of $G$ by $X$, with $\ct=|X|$.  For a set $U$
of vertices, we write $C(U)=\sum_{x\in U} C(x)$, and define $U^i=\{u\in
U\mid C(u)=i\}$.  For a list of vertices, we denote $C(x_1,\ldots,x_k)
= (C(x_1),\ldots,C(x_k))$.  We say that a graph is $r$-{\it
(semi)greedy} if every configuration of size at least $\pi(G,r)$
has a (semi)greedy $r$-solution; that is, every pebbling step in
the solution decreases (doesn't increase) the distance of the moved
pebble to $r$.  Note that any step from a cone vertex to one of its
neighbors is semigreedy.

We begin by outlining in Section \ref{WFLem} a rather new technique
for finding upper bounds on $\pi$ using weight functions.  From
there we prove pebbling number results in the case that $\ecc(r)=2$.
We prepare in Section \ref{Ecc3} preliminary lemmas that will be
used in Section \ref{PfE3} to prove pebbling results for the
$\ecc(r)=3$ case.  In Section \ref{Algo} we collect recognition
results for Pereyra and Phoenix graphs that are combined with our
pebbling number theorems to prove our main result that pebbling
numbers for split graphs can be calculated in polynomial time.  From
this analysis we learn that all split graphs with minimum degree
at least 3 are Class 0.  We end with some comments and conjectures
in Section \ref{Remarks}.

%
%
\section{The Weight Function Lemma}\label{WFLem}

In this section we describe a tool developed in \cite{HurlWFLGP}
for calculating upper bounds for pebbling numbers of graphs that
will be useful in delivering a quick proof of Theorem \ref{rootKm}.

Let $G$ be a graph and $T$ be a subtree of $G$, with at least two
vertices, rooted at vertex $r$.  For a vertex $v\in V(T)$ let $v^+$
denote the {\em parent} of $v$; i.e.  the $T$-neighbor of $v$ that
is one edge closer to $r$ (we also say that $v$ is a {\em child}
of $v^+$).  We call $T$ a {\em strategy} when we associate with it
a nonnegative, nonzero {\it weight function} $\sw:V(T)\rar\zN$ with
the property that $\sw(r)=0$ and $\sw(v^+)\ge 2\sw(v)$ for every
other vertex that is not a neighbor of $r$ (and $\sw(v)=0$ for
vertices not in $T$).  We extend $\sw$ to a function on configurations
by defining $\sw(C)=\sum_{v\in V}\sw(v)C(v)$.  Now denote by $\bT$
the configuration with $\bT(r)=0$, $\bT(v)=1$ for all $v\in
V(T)$, and $\bT(v)=0$ everywhere else.  The following was proven
in \cite{HurlWFLGP}.

\begin{lem}\label{WFL}
{\rm\bf [Weight Function Lemma]}
Let $T$ be a strategy of $G$ rooted at $r$, with associated weight
function $\sw$.  Suppose that $C$ is an $r$-unsolvable configuration
of pebbles on $V(G)$.  Then $\sw(C)\le \sw(\bT)$.
\end{lem}

The manner in which one uses this lemma to obtain a pebbling number
upper bound is as follows.  If we have several strategies
$T_1,\ldots,T_m$ of $G$, each rooted at $r$, with associated weight
functions $\sw_1,\ldots,\sw_m$ and configurations $\bT_1,\ldots,\bT_m$,
then we can define the accumulated weight function $\sw=\sum_{i=1}^m\sw_i$
and the accumulated configuration $\bT=\sum_{i=1}^m\bT_i$, and have
that $\sw(C)\le\sw(\bT)$ for every $r$-unsolvable configuration
$C$.  Moreover, if it so happens that $\sw(v)\ge 1$ for all $v\in
V-\{r\}$, then we also have $|C|\le\sw(C)$, from which follows
$\pi(G,r)\le \lfloor\sw(\bT)+1\rfloor$.

%
%
\section{Eccentricity Two}\label{Ecc2}

For a split graph $G$ define $X^r=X-\{r\}$, with $\crr=|X^r|$.

\begin{thm}\label{rootKm}
If $r \in K$ then $G$ is $r$-greedy and $\pi(G,r) = n + \crr$.
\end{thm}

\Pf
The lower bound is given by the configuration having $0$ on $r$ and
every cut vertex, $3$ on one leaf per vertex in $X^r$,
and $1$ everywhere else.  The upper bound can be proved by
using the Weight Function Lemma as follows.

For every neighbor $\rp$ of $r$ we define a strategy $T_{\rp}$.  If
$\rp\in X$ then give it weight $2$.  Include all of its neighbors
outside of $K$, giving them weight $1$ each.  If $\rp\not\in X$
then give it weight $1$.  For every vertex $s$ not yet in some
strategy (necessarily not in $K$; also $d_s \ge 2$), choose neighbors
$\sp$ and $s^\ppr$ and include $s$ in both strategies $T_{\sp}$ and
$T_{s^\ppr}$ with weight $1/2$ each.  The resulting sum of strategies
has weight $2$ on every vertex in $X^r$, and weight
$1$ everywhere else.  Hence $\pi(G,r) \le n + \crr$.

Greediness follows because every strategy used is $r$-greedy.
\pf

We recall the theorem of \cite{BlaSch,ClHoHu} that if $G$
is a $2$-connected, diameter $2$ graph then $\pi(G)=n+1$ if and
only if $G$ is a member of the following special class of graphs
$\cF$.  First, $\cF$ contains the Pyramid $P$, as well as $P-e$ for
any edge $e$ of the triangle $(a,b,c)$.  Notice that these graphs
have the following {\it separation property}: $\{a,b\}$ separates
$r$ from $c$, $\{b,c\}$ separates $q$ from $a$, and $\{a,c\}$
separates $p$ from $b$.  Next, $\cF$ is closed by adding cones over
pairs or triples from $\{a,b,c\}$.  Finally, $\cF$ is closed by
adding edges between cone vertices, provided that we maintain the
separation property.  Thus, if $G$ is a $2$-connected split graph
of diameter $2$, then $G\in\cF$ if and only if $G$ is Pereyra.  In
particular, we obtain the diameter $2$ case of Theorem \ref{rootE2},
below, when $\ct=0$ (i.e. $G$ is $2$-connected).

For a cone vertex $r$, we have two cases since $\ecc(r)\in\{2,3\}$.
We first note that, in the case $\ecc(r)=2$, every $r$-unsolvable
configuration $C$ has $C(v) \le 3$ for all $v$.  In particular, the
solution moving pebbles directly to $r$ from a vertex with $C(v)
\ge 4$ is greedy.  Recall that $\ct=|X|$.

\begin{thm}\label{rootE2}
If $r$ is a cone vertex with $\ecc(r)=2$, then $G$ is $r$-semigreedy
and $\pi(G,r) = n + \ct + \per$, where $\per=\per(G,r)$ is $1$ if
$G$ is $r$-Pereyra and $0$ otherwise.
\end{thm}

\Pf
The lower bound for non $r$-Pereyra graphs is given by the configuration
having $0$ on $r$ and every cut vertex, $3$ on one leaf per cut
vertex, and $1$ everywhere else.  For $r$-Pereyra graphs we place
$0$ on $r$, $a$, $b$, and $c$, $3$ on $p$ and $q$, and $1$ everywhere
else ($X=\mt$ because $\ecc(r)=2$).

We first prove the upper bound directly for $r$-Pereyra graphs.  If
$G$ is $r$-Pereyra then $N_r=\{a,b\}$, and since $\ecc(r)=2$ we
have $\ct=0$ and $[N_x]\cap\{a,b\}\neq\mt$ for all $x$.  If $C$ is
$r$-unsolvable of size $|C|=n+1$ then $C(r)=0$ and some $C(x)\ge
2$ with, say, $a\sim x$.  Thus $C(a)=0$, and also $C(y)\le 1$ for
all $y\in N_a$.  Now we have $n+1$ pebbles on $n-2$ vertices, which
means there must be another vertex $z$, with $b\sim z\not\sim a$,
having $C(z)\ge 2$, and so $C(b)=0$.  This puts the $n+1$ pebbles
on just $n-3$ vertices, which can only happen if $C(r,a,b,x,z) =
(0,0,0,3,3)$ and $C(y)=1$ for all other $y$.  But this allows us
to solve $r$ by moving a pebble from $x$ to $a$, from $z$ to a
common neighbor of $z$ and $a$ and then to $a$, and finally from
$a$ to $r$.  This contradiction means that every configuration of
size $n+1$ is $r$-solvable.

Next we prove the upper bound for non $r$-Pereyra $G$.  The lower
bound is given by the following two unsolvable configurations having
size $n+\ct-1$.  The first, when $r$ is the only leaf, has $0$ on
$r$ and its neighbor $\rp$, $3$ on some $x\neq\rp$, and $1$ everywhere
else.  Otherwise, the second has $0$ on $r$ and every cut vertex,
$3$ on one leaf per vertex in $X^r$, and $1$ everywhere else.

For the upper bound, as described above, the result is true for
diameter $2$ graphs, and so we may assume that $\diam(G)=3$.  This
means that $d_r\ge 2$ because, otherwise, $\ecc(r)=2$ would require
that every vertex is adjacent to the neighbor of $r$.  Moreover,
$\diam(G) = 3$ implies that there are at least two cones different
from $r$, whose neighborhoods are disjoint.  We remark first that,
with eccentricity $2$, the only nonsemigreedy move is one from
distance $1$ to distance $2$; but if a move from distance $1$ is
possible then it can move to $r$ immediately.  Therefore every
$r$-solution can be converted to one which is semigreedy.

If a cone vertex $v\neq r$ has the property that $G-v$ is $r$-Pereyra,
then we say that $v$ is {\it bad}; otherwise it is a {\it good}
cone vertex.  Notice that a bad cone vertex is necessarily a leaf
adjacent to a neighbor of $r$; in addition, it is the unique such
leaf and $d_r=2$.  Let $C$ be a configuration of size $n + \ct$.
We argue by induction (on the number of cone vertices) and contradiction
that $C$ is $r$-solvable.

Suppose that $C$ is not $r$-solvable, let $v$ be any cone vertex,
and define $\Gp=G-v$, with $\Cp=C$ on $\Gp$ and $\Cp(v)=0$.  Also
define $\ctp=\ct(\Gp)$ and $\perp=\per(\Gp)$.  Because $\Cp$ is
$r$-unsolvable on $\Gp$, we have $n-C(v)+\ct = |\Cp|<\pi(\Gp,r)$.
By induction, $\pi(\Gp,r)= (n-1)+\ctp+\perp \le (n-1)+\ct$ whether
$v$ is good or bad: if $v$ is good it holds because $\ctp\le\ct$
and $\perp=0$, and if $v$ is bad it holds because $\ctp=0$, $\perp=1$,
and $\ct=1$.  Therefore we may assume that $C(v)\ge 2$.

If $C(v) = 2$, then move a pebble from $v$ to one of its neighbors
to form $C^*$.  Then $C^*$ is a configuration on $\Gp$ of size
$n - 1 + \ct$, which by induction is $r$-solvable.  On the other
hand, if $C(v) \ge 3$, then $C(v) = 3$.  We can make the above
argument for each cone vertex; thus we may assume that $C(v) = 3$
for every cone vertex.  Hence no neighbor of $r$ is adjacent to
more than one cone vertex, and every neighbor of $r$ adjacent to
some cone vertex must have no pebble.  Furthermore, if some $x \in K$
has two pebbles then we can move pebbles greedily from $v$ to its
common neighbor $\rp$ of $r$, from $x$ to $\rp$, and then from
$\rp$ to $r$.  Hence we may assume that $C(x) \le 1$ for all $x
\in K$.

Recall that there are at least two cone vertices.  If $v$ is a cone
vertex with neighbor $\vp$ having $C(\vp) \ge 1$, then move a
pebble from another cone vertex $u$ to its common neighbor $\up$
of $r$.  Then move a second pebble from $v$ to $\vp$ to $\up$
to $r$.  Thus we must have $C(N_v) = 0$ for every cone vertex $v$.

We claim that the neighborhoods of cone vertices are pairwise
disjoint.  Indeed, suppose two cone vertices $u$ and $v$ have a
common neighbor $x$.  If there is a third cone vertex $w$ (necessarily
having $3$ pebbles), then move one to its common neighbor $\wp$
of $r$.  Then move pebbles from $u$ and $v$ to $x$, then from $x$
to $\wp$ to $r$.  Thus there are no other cone vertices.  As
mentioned above, if $u$ and $v$ are the only cone vertices then
$N_u$ and $N_v$ are disjoint.  This proves the claim.

Now we may partition $G-r$ into closed neighborhoods of cone vertices
and one extra part $\Kp$ consisting of vertices of $K$ adjacent
to no cone.  Notice that the above arguments show that $C([N_v]) =
3$ for every cone vertex $v$.  Moreover, $3 = |[N_v]| + 1$ when
$d_v = 1$ (i.e.  $\vp\in X$), and $3 \le |[N_v]|$
otherwise.  Also, $C(\Kp) \le |\Kp|$.  Hence $|C| \le n - 1 +
\ct$, a contradiction.
\pf

We finish this section with a result that will be used to prove
Theorem \ref{rootE3} below.
Define $\pi_k(G,r)$ to be the minimum number of pebbles $t$ so that
from every configuration of size $t$ one can move $k$ pebbles to $r$
(such a configuration is called {\it $k$-fold $r$-solvable}).
For example, $\pi_1(G,r)=\pi(G,r)$.

Recall that $X^r=X-\{r\}$ and $\crr=|X^r|$.
Now define $\Xrs=X-N_r-N_s$, with $\crs=|\Xrs|$.

\begin{thm}\label{2pebbK}
If $r\in K$ and $\d=\ds(G,r)$ then
$$\pi_2(G,r)=\left\{%
\begin{array}{lll} 
n + \crr + 4 & \hbox{if} & \d=1;\\ 
n + 6 -\d & \hbox{if} & 1< \d < 4;\\ 
n + 2 & \hbox{if} & \d \ge 4.\\
\end{array}%
\right.$$
\end{thm} 

\Pf
Suppose $\d=1$.
Choose $s$ to be a vertex at distance $2$ from $r$ with $d_s=\d$, 
The lower bound is given by the following configuration $C$ of size
$n+\crr+3$ that is not $2$-fold $r$-solvable:
we place $0$ pebbles on $r$ and each cut vertex, $7$ on $s$, $3$ on one leaf per 
vertex in $\Xrs$, and $1$ everywhere else.
Evidently, the only pebble that can reach $r$ comes from four that are on $s$.

For the upper bound, we assume that $C$ is a configuration of size
$n+\crr+4$ that cannot place two pebbles on $r$.
If we can place one pebble on $r$ using at most $3$ pebbling steps, then
Theorem \ref{rootKm} says we can place another on $r$ with the remaining
$n+\crr$ pebbles, so we suppose otherwise.

This means that $C(x)\le 1$ for all $x\in K$, $C(x)\le 3$ for all $x$, 
$C(N_x)=0$ for all $x\in \SP=S^2\cup S^3$, and $N_x\cap N_y=\mt$ for all
$x,y\in \SP$.
Now every $x\in \SP$ satisfies $|[N_x]|+1\ge 3\ge C(x)=C([N_x])$,
with equality if and only if $x$ is a leaf.
Hence, with $L$ denoting the set of leaves, $L^+=L\cap \SP$, 
and $U=V-\cup_{x\in \SP}[N_x]$, we have
\begin{align*}
|C| & =\sum_{x\in L^+}C([N_x])+\sum_{x\in \SP-L^+}C([N_x])+\sum_{x\in U}C(x)\\
& \le \sum_{x\in L^+}(|[N_x]|+1)+\sum_{x\in \SP-L^+}|[N_x]|+(|U|-1)\\
& \le n+\crr-1,
\end{align*}
a contradiction.

Now suppose that $1<\d<4$ --- notice that $\crr=0$ when $\d>1$.
The lower bound comes from the configuration that places $7$ on $s$,
$0$ on $r$ and $N_s$, and $1$ everywhere else, having size $n+5-\d$.
Once again, the only pebble that can reach $r$ comes from four that
are on $s$.

The very same upper bound argument above works here when $\d=2$, so
we assume that $\d=3$, whereby $C$ has size $n+3$. Suppose $C$ is
not $2$-fold $r$-solvable. Then since by Theorem \ref{rootKm} 
we have $\pi_1(G,r)=n$, it must be that:
\begin{enumerate}
\item
$C(r)=0$,
\item
$C(x)\le 1$ for every $x\in K$,
\item
if $x\in S$ and $C(x)\ge 2$ then $C(N_x)=0$,
\item\label{d:ind}
(by induction) $C(x)\ge 2$ for every $x\in S_s$, and
\item\label{e:dist}
if there exists a vertex $x\ne s$ at distance $2$ from $r$
with $d_x=\d$, then $C(s)\ge 2$.
\end{enumerate}

Now, if there exists $x\in S_s$, then by part \ref{d:ind} we have
$C(x)\ge 2$, and by part \ref{e:dist} we have $C(N_x)=0$.  Let $h\in
N_x$, $h\ne r$, and consider $\Gp=G-h$.  Notice that $\ds(\Gp,r)\ge
\d-1=2$ so that, by induction, $\pi_2(\Gp,r)= n-1+6-\ds(\Gp,r)\le
n+3=|C|$.  Thus $C$ is $2$-fold $r$-solvable, a contradiction.

Otherwise, $S_s= \emptyset$, and we can assume $K=\{r\} \cup N_s$.
It follows $n=5$, $|C|=8$, and $C$ is $2$-fold $r$-solvable,
a contraction.

Finally, suppose that $\d\ge 4$. In this case the lower bound comes
from the configuration with $3$ on $s$, $0$ on $r$, and $1$
everywhere else, having size $n+1$. Here, the only pebble that can
reach $r$ comes from two on $s$.

For the upper bound, let $C$ be a configuration of size $n+2$ that
is not $2$-fold $r$-solvable.  Since, by Theorem \ref{rootKm}, we
have $\pi(G,r)=n$, it must be that $C(r)=0$, and $C(x)\le 1$ for
every $x\in K_r$.
We will use induction on $|S|$, with the base case of $|S|=1$
(say $S=\{x\}$).
In this case we have $C(x)\ge 4$, so $C(N_x)\le 1$.
Then $C(x)\ge 8-C(N_x)$, so in either case of $C(N_x)\in\{0,1\}$
we can put two pebbles on $r$, a contradiction.
Now suppose that $|S|\ge 2$.

Let $x\in K_r$. Because $\d\ge 4$, $G-x$ is connected and has no cut
vertices different from $r$.  Denote $\dpr=\ds(G-x,r)$.  Notice that 
$\dpr \ge \d-1$ and so, by the inductive hypothesis, 
$$\pi_2(G-x,r)=\left\{%
\begin{array}{ll}
    n-1+6-3=n+2, &  \mbox{ when } \dpr=3;\\
    n-1+2=n+1, & \mbox{ when } \dpr\ge 4. \\
\end{array}%
\right. $$
This implies that if $C(x)=0$ then $C$ is $2$-fold $r$-solvable,
a contradiction.

Therefore $C(x)=1$ for every $x\in K_r$, thus $C(S)=n+2-|K_r|=|S|+3$. 
This means that in $S$ there is a vertex with at least $4$ pebbles 
or there are two vertices with at least $2$ pebbles each. In both 
cases we can place two pebbles on $r$, a contradiction which
completes the proof.
\pf

%
%
\section{Eccentricity Three}\label{Ecc3}

In the case that $\ecc(r)=3$, define 
$\dtr$ to be the set of vertices at distance $3$ from $r$, with
$\d=\ds(G,r)$, and let $s\in \dtr$ be chosen to have $d_s=\d$.
Denote by $S$ the set of cone vertices of $G$,
with $S_v=S-\{v\}$ and $\Srs=S-\{r,s\}$.
Also, let $K_v=K-N_v$, and $\Krs=K_r-N_s$.
Recall that $\Xrs=X\cap\Krs$, with $\crs=|\Xrs|$.
Now let $\Xrz$ be the set of cut vertices of $N_r$ adjacent to
some cone vertex in $S_r$, with $\crz=|\Xrz|$.
Note that $\crs>0$ implies $d_s = 1$.

Define the following four functions:
\begin{align*}
\trs(G,r) & =  n + \crs + 6 - d_r - d_s;\\
t_r(G,r) & =  n + \crs + 2 - d_r;\\
t_s(G,r) & =  n + \crs + \crz + 2 - d_s;\\
t_0(G,r) & =  n + \crs + \crz;
\end{align*}
and let $t(G,r) = \max\{t_\a(G,r)\mid \a\in\{rs,r,s,0\}\}$.
Notice that $t$ is well defined: the selection of vertex $s$ does not
change the value of $t$.
Furthermore, the choice of $S$ in the split representation of $G$
does not influence $t$ either.
Also, If $G$ is $r$-Phoenix then $d_r=2$, $\crz=\crs=0$, and $d_s=4$,
which yields $t(G,r)=n$ in this instance.

Next define the following four configurations $C_\a$ of sizes $|C_\a| = t_\a(G,r)-1$.

\begin{quote}
\begin{itemize}
\item[$\Crs$:]
$0$ on $r$, $N_r$, $N_s$, $\Xrs$, $7$ on $s$, $3$ on
one leaf per cut vertex in $\Xrs$, and $1$ everywhere else.
\item[$C_r$:]
$0$ on $r$, $N_r$, and $\Xrs$, $3$ on $s$ and on one leaf per cut
vertex in $\Xrs$, and $1$ everywhere else.
\item[$C_s$:]
$0$ on $r$, $N_s$, $\Xrs$, and $\Xrz$, $3$ on $s$ and on one leaf
per cut vertex in $\Xrs\cup\Xrz$, and $1$ everywhere else.
\item[$C_0$:]
$0$ on $r$, $\Xrs$, and $\Xrz$, $3$ on one leaf per cut vertex in
$\Xrs\cup\Xrz$, and $1$ everywhere else.
\end{itemize}
\end{quote}
Also, in the case that $G$ is $r$-Phoenix, define the configuration $C_P$ by
placing $0$ on $\{r,a,b,c\}$, $3$ on $p$ and $q$, and $1$ everywhere else.
Notice that $C_P$ witnesses that $\pi(G)\ge n+1$ for every $r$-Phoenix graph $G$.

\begin{lem}\label{unsolv}
Each $C_\a$ is $r$-unsolvable.
\end{lem}

\Pf
For $\a\in\{s,0\}$ $C_\a$ is $r$-unsolvable because the only pebbling
moves available are from the cones with $3$ pebbles to $K$, and
after those no pebbling move is available.  In $C_r$, the only move
available is from $s$ to some $v \in N_s$, and then from $v$ along
any path to some $u \in N_r$, at which point no more moves are
available.  In $\Crs$, the leaves with $3$ pebbles can only
move to their neighbors, at which point they stop.  Then $s$ can
only move $3$ to its neighbor, at which point it can travel along
any path to some neighbor of $r$ and stop there.
Finally, as mentioned in the proof of Theorem \ref{rootE2}, $C_P$ is
$r$-unsolvable on $r$-Pereyra graphs.
\pf

\begin{figure}
\centerline{\includegraphics[height=3.0in]{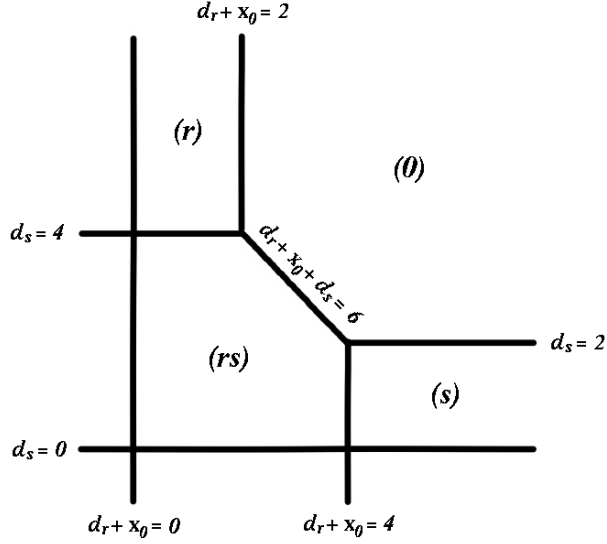}}
\caption{Graph of Cases in Lemma \ref{biggest}}\label{GraphCases}
\end{figure}

\begin{lem}\label{biggest}
With the values of $t_\a$ defined above, we list when (if and only if)
each is largest.
\begin{itemize}
\item[{\rm ($rs$)}]
$\trs \ge t_\a$ for all $\a \in \{s,r,0\}$ when
$d_s \le 4$,
$d_r + \crz \le 4$,
and
$d_r + d_s + \crz \le 6$;

\item[{\rm ($r$)}]
$t_r \ge t_\a$ for all $\a \in \{rs,s,0\}$ when
$d_s \ge 4$,
and
$d_r + \crz \le 2$;

\item[{\rm ($s$)}]
$t_s \ge t_\a$ for all $\a \in \{rs,r,0\}$ when
$d_r + \crz \ge 4$,
and
$d_s \le 2$;

\item[{\rm ($0$)}]
$t_0 \ge t_\a$ for all $\a \in \{rs,s,r\}$ when
$d_r + d_s + \crz \ge 6$,
$d_r +\crz \ge 2$,
and
$d_s \ge 2$.

\end{itemize}
\end{lem}

\Pf
Easy to check (see Figure \ref{GraphCases}).  
\pf

The next lemma shows how the function $t$ changes when some vertex
is removed. We say that a vertex $v$ has a {\em false twin} if there
exists $\vp$ non-adjacent to $v$ such that $N_v=N_{\vp}$.

\begin{lem} \label{l:removeS}
Let $v\in \Srs$.  Then

\begin{enumerate}

\item\label{parta}
If $d_v\ge 2$ then $t(G-v,r)=t(G,r)-1$.

\item\label{partb}
If $d_v= 1$ and $v$ has at least one false twin different from $r$
then $t(G-v,r)= t(G,r)-1$.

\item\label{partc} 
If $d_v= 1$ and $r$ is the only false twin vertex of $v$ then
$t(G-v,r)\le t(G,r)-1$.

\item\label{partd}
If $d_v= 1$, $v$ has no false twins, and $N_v\sse \Xrs$ then
$t(G-v,r)=t(G,r)-2$.

\item\label{parte}
If $d_v= 1$, $v$ has no false twins, and $N_v\sse \Xrz$ then
$t(G-v,r)\le t(G,r)-1$.

\end{enumerate}

\end{lem}

\Pf
This follows from Lemma \ref{biggest}.
\pf

\begin{cor} \label{c:}
If $v\in \Srs$ then $t(G-v,r)\le t(G,r)-1$.
\pf
\end{cor}

\begin{lem} \label{l:removeK}
If $d_r\ge 2$, $x \in N_r$ and $N_x\cap S=\{r\}$, then
$t(G-x,r)\le t(G,r)$.
\end{lem}

\Pf
This follows from Lemma \ref{biggest}.
\pf

\begin{lem} \label{l:becomePhoenix}
Let $G$ be non $r$-Phoenix, $\d=\ds(G,r)$, and assume there exists
$v\in S_r$ such that $\Gp=G-v$ is $r$-Phoenix. Then exactly one of
the following statements is true.
\begin{enumerate}
\item
$v$ is the only vertex of $G$ with degree $1$, and $N_v \sse N_r$.
In this case, $d_r(G)=d_r(\Gp)=2$,
$\d\ge 4$, $\crs=0$ and $\crz=1$; thus $t(G,r)=n+1$.
\item
$\d\le 3$ and $v$ is the only vertex of $\dtr$ with $d_v=\d$.
In this case, $$t(G,r)=\left\{%
\begin{array}{ll}
    n+3, & \hbox{\mbox{ if  } $\d=1$;} \\
    n+4-\d, & \hbox{\mbox{ if  } $2\le \d\le 3$.} \\
\end{array}%
\right.$$
\end{enumerate}
In both cases, if $w \neq r$ is a cone vertex of an $r$-Pyramid
of $G$ then $G-w$ is not $r$-Phoenix and $t(G-w,r)=t(G,r)-1$.
\end{lem}

\Pf
This follows from the definition of $r$-Phoenix and from Lemma \ref{biggest}.
\pf

\begin{thm}\label{rootE3}
If $r$ is a cone vertex with $\ecc(r)=3$, then $G$ is $r$-semigreedy and
$\pi(G,r) = t(G,r)+\phx(G,r)$, where $\phx(G,r)=1$ if $G$ is $r$-Phoenix and
$0$ otherwise.
\end{thm}

%
%
\section{Proof of Theorem \ref{rootE3}}\label{PfE3}

The lower bound is given by Lemma \ref{unsolv}.  The upper bound
follows by induction on $n=|V(G)|$.  The theorem is trivially
true if $n=4$.  Suppose that $G$ is a graph with at least $5$ vertices,
$r$ a cone vertex with $\ecc(r)=3$, and $C$ is a configuration on $G$
of size (without loss of generality) exactly $t = t(G,r)+\phx(G,r)$.  We
assume, for the sake of contradiction, that $C$ is not $r$-solvable;
in particular, $C(r)=0$.  The semigreediness of $G$ will follow from
moving pebbles semigreedily to a subgraph that by induction has a
resulting semigreedy solution.  Among vertices in $\dtr$,
let $s$ be chosen to have the minimum degree $\d=\ds(G,r)$ and, 
among such, having the maximum number of pebbles.

%
%
\subsection{$G$ is $r$-Phoenix}\label{rPhx}

Since $G$ is $r$-Phoenix, then $t(G,r)=n$ and so $|C|=n+1$.
Let $p \in S$ be a cone vertex of an $r$-Pyramid such that $N_p=
\{a,c\}$. It is clear that $C(p)\le 3$.
By Lemma \ref{l:removeS}(\ref{parta}), $t(G-p,r)=t(G,r)-1$. Thus, by the
inductive hypothesis, we have
$\pi(G-p,r)=t(G-p,r)+\phx(G-p,r)=t(G,r)-1+\phx(G-p,r)\le
t(G,r)=n$.

If $C(p)\le 2$, then we can move a pebble from $p$ to $N_p$,
creating a configuration $\Cp$ on $G-p$ of size $n$, which implies 
that $\Cp$, and hence $C$ is $r$-solvable, a contradiction.
So we may assume that $C(p)=3$ and, by an analogous argument, that
$C(q)=3$ where  $q$ is a cone vertex of the $r$-Pyramid such that
$N_{q}= \{b,c\}$. Moreover, we can assume that $p$ and $q$ are
the only cone vertices with degree $2$ whose neighborhoods are $\{a,c\}$
or $\{b,c\}$. It follows that the graph $G-p$ is not $r$-Phoenix,
so $\phx(G-p,r)=0$.

Then, as above, we obtain $\pi(G-p,r)=t(G,r)-1+\phx(G-p,r)=
t(G,r)-1=n-1$.  Moving a pebble from $p$ to $N_p$, we obtain a 
configuration $\Cp$ on $G-p$ of size $n+1-3+1=n-1$, which implies 
that $\Cp$, and hence $C$ is $r$-solvable, a contradiction.

%
%
\subsection{$G$ is not $r$-Phoenix}\label{notrPhx}

Thus $|C|= t(G,r)+\phx(G,r)=t(G,r)+0= t(G,r)$.

%
%
\subsubsection{$G-v$ is $r$-Phoenix for some $v \in S_r$}\label{GvrPhx}

We consider the two different cases of Lemma \ref{l:becomePhoenix}.

\begin{enumerate}
\item
The first case of Lemma \ref{l:becomePhoenix} has $t(G,r)=n+1$
and $v$ at distance $2$ of $r$; thus $C(v)\le 3$.

\begin{enumerate}
\item
If $C(v)\le 2$, we obtain a configuration $\Cp$ of $G-v$ with at
least $|C|-1= t(G,r)-1=n$ pebbles. Since $G-v$ is $r$-Phoenix,
$t(G-v,r)=n-1$, then, by the inductive hypothesis, $\pi(G-v,r)=
t(G-v,r)+\phx(G-v,r)=n-1+1=n$. This means that $\Cp$, and so $C$,
is $r$-solvable, a contradiction.

\item
If $C(v)= 3$, let $w\neq r$ be a cone vertex of an $r$-Pyramid, having
distance $2$ from $v$.
It is clear that $C(w)\le 1$; thus we obtain a configuration $|\Cp|$
of $G-w$ with at least $|C|-1= t(G,r)-1=n$ pebbles.  By the observation
at the end of Lemma \ref{l:becomePhoenix}, $t(G-w,r)=n+1-1=n$ and
$G-w$ is not $r$-Phoenix, then, by the inductive hypothesis, $\pi(G-w,r)=
t(G-w,r)+\phx(G-w,r)=n+0=n$. This means that $\Cp$, and thus $C$,
is $r$-solvable, a contradiction.
\end{enumerate}

\item
The second case of Lemma \ref{l:becomePhoenix} has two options for
$t(G,r)$, depending on the value of $\d$.

\begin{enumerate}
\item
If $d_v=\d=1$ then $|C| = t(G,r) = n+3$.  We can assume that $C(v)\le 7$.

\begin{enumerate}
\item
If $C(v)\le 6$ then since, by the inductive hypothesis, $\pi(G-v,r)
= t(G-v,r)+\phx(G-v,r) = n-1+1 = n$,  it is easy to see that $C$
is $r$-solvable, a contradiction.

\item
If $C(v)=7$ then let $w\neq r$ be a cone vertex of an $r$-Pyramid.
It is clear that $C(w)\le 1$; thus we have a configuration $\Cp$
on $G-w$ of size at least $n+2$.  By the observation at the end of
Lemma \ref{l:becomePhoenix}, $t(G-w,r)=n+3-1=n+2$.  Also, $G-w$ is
not $r$-Phoenix so, by the inductive hypothesis, $\pi(G-w,r) =
t(G-w,r)+\phx(G-w,r) = n+2+0 = n+2$.  This means that $\Cp$, and
thus $C$, is $r$-solvable, a contradiction.
\end{enumerate}

\item
If $2 \le d_v=\d \le 3$, then $|C| = t(G,r) = n+4-d$.  Let $p$ be
a cone vertex of an $r$-Pyramid such that $N_p=\{a,c\}$ with $a\in
N_r$.  Since $G-p$ is not $r$-Phoenix, by Lemma \ref{l:removeS} and
the inductive hypothesis, we can assume that $C(p)=3$.  Thus we
find the configuration $\Cp$, equal to $C$ on $G-\{w,r\}$, having
size $|C|-3 = n+4-d-3=n+1-d \ge n-2$.  By Theorem \ref{rootKm} we
have $\pi(G-\{p,r\},a)=n-2$, and so $\Cp$ is $a$-solvable, implying
that $C$ is $r$-solvable, a contradiction.
\end{enumerate}
\end{enumerate}

%
%
\subsubsection{For every $x \in S_r$, $G-x$ is not $r$-Phoenix}\label{GxNotrPhx}

Recall that $s \in \dtr$ has the maximum number of pebbles among those
vertices of $\dtr$ having $d_s=\d$.

\begin{enumerate}

\item
Some $v\in\Srs$ has $C(v)\le 2$:
We obtain a configuration $\Cp$ of $G-v$ with at
least  $|C|-1= t(G,r)-1$ pebbles. By Corollary \ref{c:},
$t(G-v,r)\le t(G,r)-1$, so by the inductive hypothesis
$\pi(G-v,r)=t(G-v,r)+\phx(G-v,r)\le t(G,r)-1+0= t(G,r)-1$. This
means that $\Cp$, and hence $C$, is solvable, a contradiction.

\item
Some $v\in\Srs$ has $C(v)\ge 4$ and every other $u\in\Srs$ has
$C(u)\ge 3$:
Notice that we can assume that $v\in \dtr$, that $C(x)\le 3$ for 
every $x\in\Srs-\{v\}$, and that $C(y)=0$ and $N_y\cap S=\{r\}$ 
for all $y\in N_r$ (in particular, $\crz=0$).
Let $\rp\in N_r$ and assume that $d_r=1$. By Theorem \ref{2pebbK} 
we have
$$\pi_2(G-r,\rp)=\left\{%
\begin{array}{lll}
n -1 + c_{\rp} + 4=n+ \crs + 4  & \hbox{if} & \d=1;\\
n -1 + 6 -\d=n+5-\d & \hbox{if} & 1< \d < 4;\\
n -1+ 2=n+1 & \hbox{if} & \d \ge 4.\\
\end{array}%
\right.$$
By Lemma \ref{biggest} (since $\crs =0$ when $\d>1$) we also have
$$t(G,r)=\left\{%
\begin{array}{lll}
n + \crs-1-1 +6=n + \crs + 4 & \hbox{if} & \d=1;\\
n-1-\d+6=n + 5-\d & \hbox{if} & 1< \d < 4;\\
n-1 +2=n+1 & \hbox{if} & \d \ge 4.\\
\end{array}%
\right.$$
Thus $C$ can place two pebbles on $\rp$, then one on $r$, a contradiction.
It follows that we can assume that $d_r\ge 2$, so that the graph $G-\rp$
is connected.

The configuration $\Cp$, the restriction of $C$ to $G-\rp$, has size
$|\Cp|=|C|=t(G,r)$. By Lemma \ref{l:removeK}, $t(G-\rp,r)\le
t(G,r)$. Since $G-\rp$ is not $r$-Phoenix, we know from the
inductive hypothesis that $\pi(G-\rp,r) =
t(G-\rp,r)+\phx(G-\rp,r)= t(G- \rp,r) \le t(G,r)$.  This
means that $\Cp$, and therefore $C$, is solvable, a contradiction.

\item
$\Srs= \mt$ or every $v \in\Srs$ has $C(v)=3$:

\begin{enumerate}

\item
$\crz\ge 1$:
Let $w$ be a leaf adjacent to $\rp\in N_r$. By Theorem \ref{rootKm},
$\pi(G-r-w,\rp)=n-2+\crs+\c$ where $\c=1$ when $d_s=1$
and $\c=0$ otherwise. We move a pebble from $w$ to $\rp$, and
consider the configuration $\Cp$, the restriction of $C$ to $G-r-w$,
of size $t(G,r)-3$.
Notice that when $d_s=1$ we have $t(G,r)-3 \ge t_s(G,r)-3 =
n+\crs+\crz-ds+2-3 = \pi(G-r-w,\rp)- \c +\crz-d_s+1\ge
\pi(G-r-w,\rp) $, and that when $d_s > 1$ we have $ t(G,r)-3 \ge
t_0(G,r)-3 = n+\crs+\crz-3 = \pi(G-r-w,\rp)- \c +\crz-1\ge
\pi(G-r-w,\rp) $. Thus, in both cases, it is possible to move
another pebble to $\rp$, a contradiction.

\item
$\crz =0$ and $\crs\ge  1$: 
Notice that in this case $C(s)\ge 3$. Let $w$ be a leaf adjacent
to $\wp\in \Krs$.

\begin{enumerate}
\item
If $w$ has no false twins, by the inductive hypothesis and by
Lemma \ref{l:removeS}(\ref{partd}), $\pi(G-w,r)=t(G-w,r)= t(G,r)-2$. We
move a pebble from $w$ to $\wp$ and consider the configuration $\Cp$,
the restriction of $C$ to $G-w$ (except with $\Cp(\wp)=C(\wp)+1$),
having size $t(G,r)-3+1= t(G,r)-2=\pi(G-w,r)$.  This makes $\Cp$,
and hence $C$, $r$-solvable, a contradiction.

\item
If $w$ has a false twin, then we can assume that $s$ has no
false twins and $C(s)=3$. Thus $w$ can be chosen as $s$ and the
proof follows as above.

\end{enumerate}

\item
$\crz =0$ and $\crs=0$: 
Recall from Lemma \ref{biggest} that in this case we have
$$t(G,r)=\left\{%
\begin{array}{llllll}
n-d_r-d_s+6, & \hbox{if} & d_r\le 4, & d_s\le 4, & d_r+d_s\le 6; &(rs) \\
n-d_r+2, & \hbox{if}& d_r\le 2, & d_s \ge 4; & &(r)  \\
n-d_s+2, & \hbox{if}& d_r\ge 4 , & d_s \le 2; &  &(s)\\
n , & \hbox{if}& d_r\ge 2, & d_s\ge 2, & d_r+d_s\ge 6. &(0) \\
\end{array}%
\right.$$
Furthermore, when $d_r=1$ we have from Theorem \ref{2pebbK} that
$|C|=t(G,r)\le \pi_2(G-r,\rp)$, where $\rp$ is the neighbor of $r$.
Thus we can place two pebbles on $\rp$ and hence solve $r$, a
contradiction.  So we will assume hereafter that $d_r\ge 2$.

\begin{enumerate}
\item
$C(N_r)> 0$:
Then there exists $\rp\in N_r$ with $C(\rp)=1$. By Theorem \ref{rootKm},
$\pi(G-r,\rp)=n-1+\c$, where $\c=1$ when $d_s=1$ and $\c=0$ otherwise. 
We consider the configuration $\Cp$, the restriction of $C$ to $G-r$
(except with $\Cp(\rp)=0$), having size $t(G,r)-1$, which is at least
$\pi(G-r,\rp)$ when $d_s>1$.  When $d_s=1$ we see that $ t(G,r)-1 \ge
t_s(G,r)-1 = n-1+2-1 =\pi(G-r,\rp)$.  In either case, $\Cp$ is $\rp$-solvable,
a contradiction.

\item
$C(N_r)=0$: 
Define the sets
\begin{align*}
\Ars &= \{x \in \Srs \mid N_x\cap N_r\neq \mt, N_x\cap N_s\neq \mt \},\\ 
A_r &= \{x \in \Srs \mid N_x\cap N_r\neq \mt, N_x\cap N_s=\mt \},\\
A_s &= \{x \in \Srs \mid N_x\cap N_r= \mt, N_x\cap N_s \neq \mt \}, {\rm\ and}\\
A_0 &= \{x \in \Srs \mid N_x\cap N_r= \mt, N_x\cap N_s = \mt \}.
\end{align*}

Of course, $K^i=\mt$ for $i\ge 4$.
Notice that, whenever $C(s)\ge 4$, $A_r \neq \mt$, $\Ars \neq \mt$, 
$K^1 \cap N_r \neq \mt$, or some pair of vertices $x,y\in \Srs$ satisfies 
$N_x \cap N_y \neq \mt$, we can assume both that $K^i=\mt$ for $i\ge 2$ and that 
either $A_0= \mt$ or the sets $[N_x]$ for $x\in A_0$ are pairwise disjoint.

We will analyze the possible intersections between the neighborhoods
of the cone vertices to compare the number of vertices and the size
of the configuration. We consider different cases depending on the
number of pebbles in $s$.  Let $\Kp=K-N(S)$.

\begin{enumerate}
\item
$6 \le C(s) \le 7 $:
In this case $A_r=\Ars=A_s= \mt$.  Thus 
$n= 1+d_r+1+d_s+ \sum_{x\in A_0} |[N_x]|+|\Kp| 
\ge 1 + d_r +1 + d_s + 3|A_0| + |K^1|$.  
We also have $C(K)=|K^1|$, and so $|C|=3|A_0|+C(s)+|K^1|$.
Then $|C|=t(G,r)\ge n-d_r-d_s + 6 \ge 1 + d_r +1 + d_s +3|A_0| +|K^1|
-d_r-d_s + 6 = |C|-C(s)+8$.  Thus $C(s) \ge 8$, a contradiction.

\item
$4 \le C(s) \le 5$:
In this case $A_r=\Ars= \mt$.
Moreover, $K^1\sse N_s$, $|A_s|+|K^1|\le 1$, and $N_x\cap N_y=\mt$
for all $\{x,y\}\sse \Srs$ ($x\neq y$).
This means that $|C|=3|A_0|+3|A_s|+|K^1|+C(s)$ and
$n\ge 1+d_r+1+d_s+\sum_{x\in A_0}|[N_x]|+|A_s|
\ge d_r+d_s+3|A_0|+2+|A_s|$.
Together these imply that
$|C|=t(G,r)\ge n-d_r-d_s+6 \ge 3|A_0|+8+|A_s|
= |C|-2|A_s|-|K^1|-C(s)+8$,
and hence $C(s)\ge 8-2|A_s|-|K^1|\ge 6$, a contradiction.

\item
$2 \le C(s) \le 3$:

\begin{itemize}
\item[I.]
If $A_{r}\neq \mt$: 
Then $\Ars=A_{s} = \mt$, $K^i=\mt$ for $i\ge 2$,
and $K^1\sse N_{A_r}-N_r-N_s$.

\begin{itemize}
\item[\str]
If $|A_r|\le 2$ then
$n \ge 1+d_r+1+d_s +\sum_{x\in A_0}|[N_x]|+|A_r|+|K^1|
\ge d_r+d_s+2+3|A_0|+|A_r|+|K^1|$.
Also
$3|A_0|+3|A_r|+C(s)+|K^1| = |C|
\ge n-d_r-d_s+6
\ge 8+3|A_0|+|A_r|+|K^1|$,
which implies the contradiction that
$C(s)\ge 8-2|A_r|\ge 4$.

\item[\sstr]
If $|A_r|\ge 3$ then $K^1=\mt$ and
$n \ge 1+1+d_s +\sum_{x\in A_0\cup A_r}|[N_x]|
\ge d_s+2+3|A_0|+3|A_r|$.
Thus
$3|A_0|+3|A_r|+C(s) = |C|
\ge n-d_s+2
\ge 4+3|A_0|+3|A_r|$,
which implies the contradiction that $C(s)\ge 4$.
\end{itemize}

\item[II.]
If $A_{r}= \mt$ and $\Ars \neq \mt$: 
Then $\Ars$ contains
exactly one vertex $w$ and $K^1 \sse N_w$.  In this case we
see that the sets $[N_r]$, $[N_x]$ (for all $x\in A_0$), $K^1$
and $[N_s]$ are pairwise disjoint. Thus $|C| \le  3+ |K^1| +
3|A_0| + C(s)$ and $|C|=t(G,r)\ge n-d_r-d_s + 6 \ge 1+ d_r+ 1+
|K^1|+ 3|A_0| +1+ d_s -d_r-d_s +6$, which implies $C(s)\ge 6$, 
a contradiction.

\item[III.]
If $A_{r}= \Ars=\mt$:
Let $\rp\in N_r$ and consider $\Gp=G- ([N_r]- \rp)$.
Notice that if $\d=1$ then $c_{\rp}=\crs+1$, with $c_{\rp}=c_{rs}=0$
otherwise.
By Theorem \ref{2pebbK},
$$\pi_2(\Gp,\rp)=\left\{%
\begin{array}{lll}
n- d_r+ \crs + 5  & \hbox{if} & \d=1;\\
n -d_r + 6 -\d & \hbox{if} & 1< \d < 4\\
n -d_r+ 2 & \hbox{if} & \d \ge 4.\\
\end{array}%
\right.$$
Since $C(N_r)=0$, the restriction of $C$ to $\Gp$ has size
$$t(G,r)=\left\{%
\begin{array}{lll}
n + \crs + 5 -d_r & \hbox{if} & \d=1, d_r\le 4;\\
n + \crs +1 & \hbox{if} & \d=1, d_r\ge 4;\\
n+4-d_r & \hbox{if} & \d=2 , d_r\le 4;\\
n & \hbox{if} & \d=2 , d_r\ge 4;\\
n+3-d_r & \hbox{if} & \d=3 , d_r\le 3;\\
n & \hbox{if} & \d=3 , d_r\ge 3;\\
n+1 & \hbox{if} & \d \ge 4, d_r =1;\\
n & \hbox{if} & \d \ge 4, d_r\ge 2.\\
\end{array}%
\right.$$
Thus $C$ is $2$-fold $\rp$-solvable, hence $r$-solvable, a contradiction.

\end{itemize}

\item
$C(s)\le 1$: 
In this case, we have a configuration $\Cp$ (the restriction of 
$C$ to $G-s$) of size at least $|C|-1 =t(G,r)-1$ on the graph $G-s$. 
We will show that $\pi (G-s,r)\le t(G,r)-1$, implying that $\Cp$,
and hence $C$ is $r$-solvable, a contradiction.

\begin{itemize}
\item[I.]
If $r$ has eccentricity $2$ in $G-s$ and $G-s$ is not Pereyra:
Then $\pi(G-s,r)=n-1$. On the other hand $t(G,r)\ge n$.

\item[II.]
If $r$ has eccentricity $2$ in $G-s$ and $G-s$ is $r$-Pereyra: 
Then $\pi(G-s,r)=n-1+1=n$ and $d_r=2$.  
Furthermore, $d_s \le 3$ because $G$ is not $r$-Phoenix.
Hence $t(G,r)=n-2-d_s+6\ge n +1$.

\item[III.]
If $r$ has eccentricity $3$ in $G-s$:
Then, by the inductive
hypothesis, $\pi (G-s,r)=t(G-s,r)$, since we know that $G-s$
is not $r$-Phoenix. Let $\dpr=\ds(G-s,r)$ and notice that, since any cone
vertex of $\Srs$ has $3$ pebbles and $s$ has just one pebble,
then $d_s < \dpr$. We have from Lemma \ref{biggest} that
$$t(G-s,r)=\left\{%
\begin{array}{lllll}
n-d_r-\dpr+5 & \hbox{if} & d_r\le 4,\ \dpr\le 4, & &(rs)^\pr \\
 & & {\rm and\ } d_r+\dpr\le 6;\\
n-d_r+1 & \hbox{if}& d_r\le 2,\ \dpr \ge 4; & &(r)^\pr  \\
n-\dpr+1 & \hbox{if}& d_r\ge 4,\ \dpr \le 2; &  &(s)^\pr\\
n-1  & \hbox{if}& d_r\ge 2,\ \dpr\ge 2, & &(0)^\pr \\
 & & {\rm and\ } d_r+\dpr\ge 6.
\end{array}%
\right.$$

Observe that the only possible change of cases from $G$ to $G-s$
is from $(rs)$ to $(r)^\pr$ or $(0)^\pr$, or from $(s)$ to $(0)^\pr$.  
It is easy to see that in all cases, $t(G-s,r)\le t(G,r)-1$.  
\end{itemize}

\end{enumerate}
\end{enumerate}
\end{enumerate} 
\end{enumerate}

This completes the proof.
\pf

For $n = 2m$ ($+1$ if $n$ is odd), define the {\it sun} $S_n$, to
be the split graph with $|K|=m$ and $m$ leaves matched with the
vertices of $K$ (and an extra leaf joined to $K$ if necessary).
According to Theorem \ref{rootE3} we have $\pi(S_n) = n + (m - 2)
+ (6 - 1 - 1) = \lf 3n/2\rf + 2$, showing that the pebbling
bound for diameter $3$ graphs given in \cite{PoStYe} is tight.

%
%
\section{Algorithms}\label{Algo}

We begin with a key construction for finding a Pyramid in a split
graph $G$.  Suppose that $r$ is a cone vertex of $G$ with $d_r=2$.
Then let $X$ be the set of cut vertices of $G$, $W$ be the set of
degree $2$ vertices of $G$ whose neighbors are in $G-X$ and define
the graph $H=H(G)$ to have vertices $\cup_{v\in W}N_v$ and edges
$\{N_v\}_{v\in W}$.

\begin{thm}\label{rpereyra}
Given a split graph $G$ and root $r$, recognizing if $G$ is $r$-Pereyra 
can be done in linear time.
\end{thm}

\Pf
Of course $G$ being $r$-Pereyra requires $d_r=2$.  The graph $H=H(G)$
takes linear time to construct.  Then $G$ is $r$-Pereyra if and
only if $H$ has a triangle including the edge $N_r$, which can be
checked in linear time.
\pf

\begin{cor}\label{linear}
Calculating $\pi(G,r)$ when $G$ is a split graph with root $r$
can be done in linear time.
\end{cor}

\Pf
The set of cut vertices $X$ of $G$ is the neighborhood of the degree
$1$ cone vertices, and so can be calculated in linear time at the
start.  For $r \in K$, Theorem \ref{rootKm} determines $\pi(G,r)$
immediately.  For a cone vertex $r$, we calculate its eccentricity
in linear time via breadth-first search.  If its eccentricity is
$2$ then Theorem \ref{rootE2} determines $\pi(G,r)$ in linear time
from recognizing if it is $r$-Pereyra or not.  Otherwise, we have
$\ecc(r)=3$.  In the breadth-first search we also learned of all
cone vertices $s$ at distance $3$ from $r$.  As we encounter each
such $s$ we keep track of the one having least degree.  At the end
we calculate $t(G,r)$ immediately from Lemma \ref{biggest} and find
$\pi(G,r)$ via Theorem \ref{rootE3}.
\pf

Finding a triangle in a graph is a well-known problem in combinatorial
optimization.  The best known algorithm is found in \cite{AlYuZw},
below.  Let $\w\cong 2.376$ be the exponent of matrix multiplication,
and define $\b=2\w/(\w+1)\cong 1.41$.

\begin{alg}\label{triangles}
{\rm [\cite{AlYuZw}, Theorem 3.5]}
Deciding whether a graph $G$ with $m$ edges contains a triangle, 
and finding one if it does, can be done in $O(m^\b)$ time.
\end{alg}

\begin{thm}\label{pereyra}
Given a split graph $G$, recognizing if it is Pereyra can be
done in $O(n^{1.41})$ time.
\end{thm}

\Pf
We define $H=H(G)$ as above and see that $G$ is Pereyra if and only
if $H$ has a triangle.  Then Algorithm \ref{triangles} decides this
in $O(n^{1.41})$ time, since the number of edges of $H$ is at most
the number of vertices of $G$.
\pf

\begin{thm}\label{SplitPi}
If $G$ is a diamter 3 split graph then $\pi(G)$ is given as follows.
\begin{enumerate}
\item 
If $\ct\ge 2$ then
\begin{quote}
$\pi(G)=n+\ct+2$.  
\end{quote}
\item 
If $\ct=1$ then
\begin{quote}
$\pi(G)=\left\{%
\begin{array}{ll}
n+5-\ds & \hbox{if\ } r \hbox{\ is a leaf with\ } \ecc(r)=3
        \hbox{\ and\ } \ds=\ds(G,r)\le 4;\\
n+1 & \hbox{otherwise.}
\end{array}%
\right.$
\end{quote}
\item 
If $\ct=0$ then
\begin{quote}
$\pi(G)=\left\{%
\begin{array}{ll}
n+4-\ds & \hbox{if there is a cone vertex\ } r \hbox{\ with\ }d_r=2, \ecc(r)=3\\
        & \quad \hbox{and\ } \ds=\ds(G,r)\le 3;\\
n+1 & \hbox{if no such\ } r \hbox{\ exists and\ } G \hbox{\ is Pereyra};\\
n & \hbox{otherwise.}
\end{array}%
\right.$
\end{quote}
\end{enumerate}
\end{thm}

\Pf
First we remark that $\crr\le \ct$.  Hence we know that
$\pi(G)=\pi(G,r)$ for some cone root $r$.

If $\ct\ge 2$ then there exist leaves $r$ and $s$ at distance $3$
from each other (in fact, if $r$ is a leaf then so is $s$).  For
every such $r$ and $s$ we have $t(G,r)=\trs(G,r)=n+\crs+6-d_r-d_s$
from Lemma \ref{biggest}.  Also, $\crs=\ct-2$ and $d_r=d_s=1$, so
that $t(G,r)= n+\ct+2$ when $r$ is a leaf.  When $\ecc(r)=3$ but
$r$ is not a leaf, we see that $t(G,r)\le n+\ct+2$ (with equality
if and only if $d_r=2$, $d_s=1$, and $x_0=0$).  Finally, when
$\ecc(r)=2$ we have from Theorem \ref{rootE2} that
$\pi(G,r)=n+\ct+\per<n+\ct+2$.  Hence Theorem \ref{rootE3} implies
$\pi(G)=n+\ct+2$.

If $\ct=1$ then $G$ is not Phoenix.  When $\ecc(r)=2$, $G$ is not
Pereyra, and so Theorem \ref{rootE2} gives $\pi(G,r)=n+1$.  When
$\ecc(r)=3$, the cut vertex $v$ is a neighbor of either $r$ or $s$,
and so $\crs=0$.  The function $t_r=n+2-d_r$ is maximized at $n+1$
when $r$ is a leaf, so $\pi(G)\ge n+1$.  Obviously $t_0\le n+1$,
and $t_s=n+\crz+2-d_s\le n+1$, since $d_s=1$ implies $\crz=0$.  The
function $\trs$ is also maximized when $r$ is a leaf.  Indeed, if
$v\not\in N_r$ then $s$ is a leaf.  Then with $\rp=s$ having
corresponding $\sp\in D_3(\rp)$ we have $\trsp\ge \trs$ because
$d_\rp=d_s=1$ and $d_\sp\le d_r$.  So we may assume that $v\in N_r$.
If $r$ is not a leaf then let $w$ be a leaf.  But then with $\rp=w$
having corresponding $\sp\in D_3(\rp)$ we have $\trsp>\trs$ since
$d_\rp<d_r$ and $d_\sp=d_s$.  Thus we have $\pi(G) \ge n+5-d_s$
when $r$ is a leaf and $s\in D_3(r)$ with $d_s=\ds$.

Finally, if $\ct=0$, we note from Lemma \ref{biggest} and Theorem
\ref{rootE3} that the only way to have $\pi(G,r)\ge n+1$ when some
cone vertex $r$ has $\ecc(r)=3$ is either via $\trs$ (with $d_r=2$
and $d_s\le 3$) or if $G$ is $r$-Phoenix.  When a cone vertex $r$
has $\ecc(r)=2$ then we have $\pi(G,r)=n+1$ if $G$ is $r$-Pereyra,
by Theorem \ref{rootE2}.  Thus $\pi(G,r)=n$ in all other cases.

The above description can be reorganized as follows.  Suppose that
there is no cone vertex $r$ with $d_r=2$ and $s\in D_3(r)$ with
$d_s=\ds(G,r)\le 3$.  If $G$ is Pereyra then it is $r$-Pereyra for
some cone vertex $r$ with $d_r=2$.  Now we know that either $\ecc(r)=2$
or $\ds(G,r)\ge 4$, the latter case of which makes $G$ $r$-Phoenix.
In either case we get $\pi(G,r)=n+1$.
\pf

\begin{cor}\label{linpi}
Calculating $\pi(G)$ when $G$ is a split graph can be done in
$O(n^{1.41})$ time.
\end{cor}

\Pf
Recall that we discover the value of $\ct$ in linear time.  So if
$\ct\ge 2$ then $\pi(G)=n+\ct+2$.  When $\ct=1$ we let $r$ be any
leaf of $G$.  Using breadth-first search from $r$ we discover if
$D_3(r)\neq\mt$ and, if so, find $s\in D_3(r)$ with $d_s=\ds(G,r)$.
Thus, in linear time we know $\pi(G)$.

Now, if $\ct=0$, we describe a linear algorithm either to find
a cone vertex $r$ with $d_r=2$ and some $s\in D_3(r)$ having $d_s\le
3$ or to conclude that none exist.

For ease of notation, we write $d_i$ for $d_{v_i}$ and $N_i$ for
$N_{v_i}$.  In linear time we can reorder the vertices of $G$ so
that $d_i=2$ for $1\le i\le k$ and $d_i=3$ for $k+1\le i\le l$.
Initialize $\l(i)$ to be empty for every vertex $v_i$ of $G$.  Then
for each $i\le k$ we add $i$ to $\l(j)$ for each $v_j\in N_i$ and
check the size of $L_i=\cup_{v_j\in N_i}\l(j)$.  If $|L_i|<i$ then
there is some $j<i$ such that $N_i\cap N_j=\mt$ --- choose any $j\in
\{1,\ldots,i\}-L_i$.  In this case we set $r=v_j$ and $s=v_i$ and
quit; otherwise we continue.  Then for each $k+1\le i\le l$ we only
check the size of $L_i=\cup_{v_j\in N_i}\l(j)$.  If $|L_i|<k$ then
there is some $j\le k$ such that $N_i\cap N_j=\mt$ --- choose any
$j\in \{1,\ldots,k\}-L_i$.  In this case we set $r=v_j$ and $s=v_i$
and quit; otherwise we continue.  If we have not found $r$ and $s$
by now, they do not exist.  This algorithm is linear because of the
bounded degrees.

If $r$ and $s$ were found then $\pi(G)=n+6-d_r-d_s$.  If no such
$r$ and $s$ exist, we use Theorem \ref{pereyra} to discover if $G$
is Pereyra, which takes  $O(n^{1.41})$ time.  If it is then
$\pi(G)=n+1$, otherwise $\pi(G)=n$.
\pf

%
%
\section{Remarks}\label{Remarks}

We begin by noting the following corollary to Theorem \ref{SplitPi}.
Let $\k(G)$ denote the connectivity of $G$.

\begin{cor}\label{3Conn}
If $G$ is a split graph with $\d(G)\ge 3$ then $G$ is Class 0.
\end{cor}

\Pf
The first two instances of the $\ct=0$ case of Theorem \ref{SplitPi}
require $\d(G)=2$.
\pf 

Note that this implies that every $3$-connected split graph is Class 0.  
The analogous result with ``split'' replaced by ``diameter two''
was proven in \cite{ClHoHu}.  The full characterization of diameter
two, connectivity two, non Class 0 graphs in \cite{ClHoHu} involves
the appearance of a Pyramid, whereas for connectivity two, non Class
0 split graphs, Pereyra and Phoenix graphs play a significant role.

With the similarities in structure and function mentioned above
between Pyramid and Pereyra graphs, one wonders two things.  First,
in the diameter $2$ case, it is possible to add edges between twin
cone vertices (thus leaving the class of split graphs) without
changing the pebbling number; is the same true for diameter $3$?
Second, might there be a family of graphs that plays for diameter
$4$ graphs the same role played by Pyramid and Pereyra graphs for
diameters $2$ and $3$, at least in the case that the root $r$ has 
eccentricity $2$?

It is interesting that, while one can calculate the pebbling number
of a diameter two graph in polynomial time, it was shown in
\cite{CuLeSiTa} that it is {\sf NP}-complete to decide if a given
configuration on a diameter two graph can solve a fixed root.  (The
same was proven more recently for planar graphs in \cite{CuDiLe}
--- the problem is polynomial for planar diameter two graphs.)  In
that context we offer the following.

\begin{prb}
Let $C$ be a configuration on a split graph $G$ with root $r$.
Is it possible in polynomial time to determine if $C$ is $r$-solvable?
\end{prb}

We also offer the following two conjectures.

\begin{cnj}
If $G$ is chordal then $\pi(G)$ can be calculated in polynomial time.
\end{cnj}

\begin{cnj}
For fixed $d$, if $\diam(G)=d$ then $\pi(G)$ can be calculated in polynomial time.
\end{cnj}

At the very least we believe that, for a chordal or fixed diameter graph $G$,
it can be decided in polynomial time whether or not $G$ is Class $0$.

%
%
\section{Acknowledgements}\label{Ack}

The third author thanks the Fulbright International Exchange Program
for their support and assistance, the National University of La
Plata for their hospitality, the Ferrocarril General Roca for a
surprise trip to Pereyra, and his coauthors for their generosidad
extraordinaria y mate delicioso.

\bigskip

%
%

\end{document}